# Model updating using sum of squares (SOS) optimization to minimize modal dynamic residuals

[1] Dan Li, [1] Xinjun Dong, [1,2] * Yang Wang

[1] School of Civil and Environmental Engineering, Georgia Institute of Technology, Atlanta, GA, USA
[2] School of Electrical and Computing Engineering, Georgia Institute of Technology, Atlanta, GA, USA
*yang.wang@ce.gatech.edu

Abstract: This research studies structural model updating through sum of squares (SOS) optimization to minimize modal dynamic residuals. In the past few decades, many model updating algorithms have been studied to improve the similitude between a numerical model and the as-built structure. Structural model updating usually requires solving nonconvex optimization problems, while most off-the-shelf optimization solvers can only find local optima. To improve the model updating performance, this paper proposes the SOS global optimization method for minimizing modal dynamic residuals of the generalized eigenvalue equations in structural dynamics. The proposed method is validated through both numerical simulation and experimental study of a four-story shear frame structure.

## 1   Introduction

During the past few decades, numerical simulation has been widely applied to the modeling and analysis of various structures, including concrete structures [1], steel structures [2], and composite structures [3]. However, due to the complexity and large scale of civil structures, numerical models can only simulate the as-built structure with limited accuracy. In particular, material properties and boundary conditions of the as-built structure usually differ from their representations in a corresponding numerical model [4, 5]. To address this limitation, structural model updating can be performed to improve the similitude between an as-built structure and its numerical model. Meanwhile, structural model updating has been adopted for structural health monitoring in detecting the damage of structural components [6].

Various model updating methods and algorithms have been developed and applied in practice. Many of these methods utilize the modal analysis results from field testing data [7]. Structural parameter values of a numerical model are updated by forming an optimization objective function that minimizes the difference between experimental and simulated results. For example, early researchers utilized the experimentally measured eigenfrequencies, attempting to fine-tune the simulation model parameters so that the model provides similar eigenfrequencies. However, for predicting the simultaneous response at various locations of a structure with multiple degrees of freedom, it was later revealed that only the eigenfrequency data is not sufficient [8]. To this end, the modal dynamic residual approach achieves structural model updating by forming an optimization problem that minimizes the residuals of the generalized eigenvalue equations in structural dynamics [9-11]. Nevertheless, despite past efforts, these optimization problems in structural model updating are generally nonconvex. Most off-the-shelf optimization algorithms can only find some local optima, while providing no knowledge on the global optimality. In pursuing a better solution to these nonconvex problems, researchers either use randomized multiple starting points for the search [12-14] or resort to stochastic searching methods [15-19]. Nevertheless, these methods can only improve the probability of finding the global optimum; none of them can guarantee to find the global optimum.

Although the optimization problem in structural model updating is generally nonconvex, the objective function, as well as equality and inequality constraints, are usually formulated as polynomial functions. This property enables the possibility of finding the global optimum of the nonconvex problem by sum of squares (SOS) optimization method. The SOS method tackles the problem by decomposing the original objective function into SOS polynomials to find the best lower bound of the objective function. This makes the problem more solvable. Using $\mathbb{N}$ to represent nonnegative integers, an SOS polynomial $s(\mathbf{x})$ of $\mathbf{x} \in \mathbb{R}^n$ with degree of $2t$, $t \in \mathbb{N}$, represents a polynomial that can be written as the sum of squared polynomials:

$$s(\mathbf{x}) = \mathbf{z}^T \mathbf{Q} \mathbf{z} = \mathbf{z}^T \mathbf{L} \mathbf{L}^T \mathbf{z} = \sum_j \left(\mathbf{L}_j^T \mathbf{z}\right)^2, \qquad \mathbf{Q} \succcurlyeq 0 \tag{1}$$

where $n_\mathbf{z} = \binom{n+t}{n}$ is the number of $n$-combinations from a set of $n+t$ elements [20], $\mathbf{z} = (1, x_1, x_2, \cdots, x_n, x_1^2, x_1 x_2, \cdots, x_{n-1} x_n^{t-1}, x_n^t)^T \in \mathbb{R}^{n_\mathbf{z}}$ represents all the base monomials of degree less than or equal to $t$; $\mathbf{Q} \in \mathbb{R}^{n_\mathbf{z} \times n_\mathbf{z}}$ is a positive semidefinite matrix (denoted as $\mathbf{Q} \succcurlyeq 0$). Through many decomposition methods, such as eigenvalue decomposition and Cholesky decomposition, the positive semidefinite matrix $\mathbf{Q}$ can be decomposed as $\mathbf{L}\mathbf{L}^T$ with $\mathbf{L} \in \mathbb{R}^{n_\mathbf{z} \times n_\mathbf{z}}$. The $j$-th column of matrix $\mathbf{L}$ is denoted as $\mathbf{L}_j$ in Eq. (1).

In recent years, researchers in mathematical communities have applied the SOS method to calculate the global bounds for polynomial functions [21, 22]. It has also been reported that the dual problem of the SOS optimization formulation provides information about the minimizer of the original polynomial function [23-25]. Utilizing the primal and dual problems of SOS optimization, we found that for nonconvex model updating problems using the modal dynamic residual formulation, the global optimum can be reliably solved. This paper reports the findings.

The rest of the paper is organized as follows. Section 2 presents the formulation of the modal dynamic residual approach for model updating. Section 3 describes the SOS optimization method and its application on modal dynamic residual approach. Section 4 shows numerical simulation and a laboratory experiment validating the proposed SOS method. In the end, Section 5 provides a summary and discussion.

## 2   Modal dynamic residual approach

The objective of structural model updating is to identify accurate physical parameter values of an as-built structure. For brevity, we only provide formulation that updates stiffness values (although the formulation can be easily extended for updating mass and damping). Consider a linear structure with $N$ degrees of freedom (DOFs). The stiffness parameter updating is represented by a vector variable $\boldsymbol{\theta} \in \mathbb{R}^{n_\boldsymbol{\theta}}$, where each entry $\theta_i$ is the relative change from the initial/nominal value of a selected stiffness parameter (to be updated). The overall stiffness matrix can be written as an affine matrix function of the updating variable $\boldsymbol{\theta}$:

$$\mathbf{K}(\boldsymbol{\theta}) = \mathbf{K}_0 + \sum_{i=1}^{n_\boldsymbol{\theta}} \theta_i \mathbf{K}_{0,i} \tag{2}$$

where $\mathbf{K}_0 \in \mathbb{R}^{N \times N}$ denotes the initial stiffness matrix; $\mathbf{K}_{0,i} \in \mathbb{R}^{N \times N}$ denotes the $i$-th (constant) stiffness influence matrix corresponding to the updating variable $\theta_i$. Finally, $\mathbf{K}(\boldsymbol{\theta}): \mathbb{R}^{n_\boldsymbol{\theta}} \to \mathbb{R}^{N \times N}$ represents that the structural stiffness matrix is written as an affine matrix function of vector variable $\boldsymbol{\theta} \in \mathbb{R}^{n_\boldsymbol{\theta}}$. Because some stiffness parameters may not need updating, it is not required that $\mathbf{K}_0 = \sum_{i=1}^{n_\boldsymbol{\theta}} \mathbf{K}_{0,i}$.

Modal dynamic residual approach is adopted here to accomplish model updating [9-11]. The approach attempts to minimize the residuals of the generalized eigenvalue equations in structural dynamics. The residuals are calculated using matrices generated by the numerical model in combination with experimentally-obtained modal properties. Obtained through dynamic modal testing, such experimental results usually include the first few resonance frequencies ($\omega_i, i = 1, \cdots, n_{\text{modes}}$) and corresponding mode shapes. Here $n_{\text{modes}}$ denotes the number of experimental modes. For mode shapes, the experimental results can only contain entries that are associated with the DOFs instrumented with sensors. These experimentally obtained mode shape entries are grouped as $\boldsymbol{\psi}_{i,\text{m}} \in \mathbb{R}^{n_{\text{instr}}}$ where $n_{\text{instr}}$ denotes the number of instrumented DOFs, with the maximum magnitude normalized to be 1. The entries corresponding to the unmeasured DOFs, $\boldsymbol{\psi}_{i,\text{u}} \in \mathbb{R}^{N-n_{\text{instr}}}$, are unknown and need to be treated as optimization variables. The optimization problem of the modal dynamic residual approach is formulated as follows, with updating variables $\boldsymbol{\theta} \in \mathbb{R}^{n_\theta}$ and unmeasured mode shape entries $\boldsymbol{\psi}_\text{u} = (\boldsymbol{\psi}_{1,\text{u}}, \boldsymbol{\psi}_{2,\text{u}}, \cdots, \boldsymbol{\psi}_{n_{\text{modes}},\text{u}})^\text{T} \in \mathbb{R}^{(N-n_{\text{instr}}) \cdot n_{\text{modes}}}$ as the optimization variables:

$$\begin{aligned}
\underset{\boldsymbol{\theta}, \boldsymbol{\psi}_\text{u}}{\text{minimize}} \quad & \sum_{i=1}^{n_{\text{modes}}} \left\| [\mathbf{K}(\boldsymbol{\theta}) - \omega_i^2 \mathbf{M}] \begin{Bmatrix} \boldsymbol{\psi}_{i,\text{m}} \\ \boldsymbol{\psi}_{i,\text{u}} \end{Bmatrix} \right\|_2^2 \\
\text{subject to} \quad & \mathbf{L}_\theta \leq \boldsymbol{\theta} \leq \mathbf{U}_\theta \\
& \mathbf{L}_{\boldsymbol{\psi}_\text{u}} \leq \boldsymbol{\psi}_\text{u} \leq \mathbf{U}_{\boldsymbol{\psi}_\text{u}}
\end{aligned} \quad (3)$$

where $\|\cdot\|_2$ denotes the $\mathcal{L}_2$-norm; $\mathbf{M}$ is the mass matrix, which is considered accurate. It's implied that both $\mathbf{K}(\boldsymbol{\theta})$ and $\mathbf{M}$ are reordered by the instrumented and un-instrumented DOFs in $\boldsymbol{\psi}_{i,\text{m}}$ and $\boldsymbol{\psi}_{i,\text{u}}$. Constants $\mathbf{L}_\theta$ and $\mathbf{L}_{\boldsymbol{\psi}_\text{u}}$ denote the lower bounds for vectors $\boldsymbol{\theta}$ and $\boldsymbol{\psi}_\text{u}$, respectively; constants $\mathbf{U}_\theta$ and $\mathbf{U}_{\boldsymbol{\psi}_\text{u}}$ denote the upper bounds. Note that the sign "≤" is overloaded to represent entry-wise inequality.

Although the box constraints in Eq. (3) define a convex feasible set, the objective function is a fourth order polynomial which is nonconvex in general. A special case when all DOFs are measured, i.e. where $\boldsymbol{\psi}_\text{u}$ vanishes, leads to a convex objective function, because inside the $\mathcal{L}_2$-norm is an affine function of the updating variables $\boldsymbol{\theta}$. According to composition rule [26], this objective function without $\boldsymbol{\psi}_\text{u}$ is convex. However, in practice, usually not all DOFs are instrumented/measured, i.e. $\boldsymbol{\psi}_\text{u}$ exists in the objective function, rendering nonconvexity. When the problem is nonconvex, off-the-shelf optimization algorithms can only find some local optima, without guarantee of global optimality. To address the challenge, we propose the SOS optimization method that can recast the problem in Eq. (3) into a convex optimization problem and thus, reliably find the global optimum.

## 3 SOS optimization method

### 3.1 Primal problem

The sum of squares (SOS) optimization method is applicable to polynomial optimization problems. To represent a polynomial function of a vector variable $\mathbf{x} \in \mathbb{R}^n$, we use $\boldsymbol{\alpha} = (\alpha_1, \alpha_2, \cdots, \alpha_n) \in \mathbb{N}^n$ to denote the corresponding nonnegative integer-valued powers. In addition, we use a compact notation $\mathbf{x}^{\boldsymbol{\alpha}} = x_1^{\alpha_1} x_2^{\alpha_2} \cdots x_n^{\alpha_n}$ to represent the corresponding base monomial. The degree of a base monomial $\mathbf{x}^{\boldsymbol{\alpha}}$ equals $\sum_{i=1}^n \alpha_i$. With $c_{\boldsymbol{\alpha}} \in \mathbb{R}$ as the real-valued coefficient, a polynomial $f(\mathbf{x}): \mathbb{R}^n \to \mathbb{R}$ is defined as a linear combination of monomials:

$$f(\mathbf{x}) = \sum_{\boldsymbol{\alpha}} c_{\boldsymbol{\alpha}} \mathbf{x}^{\boldsymbol{\alpha}} = \sum_{\boldsymbol{\alpha}} c_{\boldsymbol{\alpha}} x_1^{\alpha_1} x_2^{\alpha_2} \cdots x_n^{\alpha_n} \quad (4)$$

Note that the summation is indexed on different power vectors $\boldsymbol{\alpha}$ that define different base monomials $x_1^{\alpha_1} x_2^{\alpha_2} \cdots x_n^{\alpha_n}$. Based on these notations, we now consider a general optimization problem with a polynomial objective function $f(\mathbf{x})$ and multiple polynomial inequalities $g_i(\mathbf{x}) \geq 0, i = 1,2,\cdots l$:

$$\begin{aligned} \underset{\mathbf{x}}{\text{minimize}} \quad & f(\mathbf{x}) = \sum_{\boldsymbol{\alpha}} c_{\boldsymbol{\alpha}} \mathbf{x}^{\boldsymbol{\alpha}} \\ \text{subject to} \quad & g_i(\mathbf{x}) = \sum_{\boldsymbol{\beta}_i} h_{\boldsymbol{\beta}_i} \mathbf{x}^{\boldsymbol{\beta}_i} \geq 0, \quad (i = 1,2,\cdots l) \end{aligned} \tag{5}$$

where $f(\mathbf{x}): \mathbb{R}^n \to \mathbb{R}$ and $g_i(\mathbf{x}): \mathbb{R}^n \to \mathbb{R}$ are polynomials with degree $d$ and $e_i \in \mathbb{N}$, respectively. Similar to $\mathbf{x}^{\boldsymbol{\alpha}}$ and $c_{\boldsymbol{\alpha}}$ in $f(\mathbf{x})$, $\mathbf{x}^{\boldsymbol{\beta}_i}$ represents a base monomial with the power of $\boldsymbol{\beta}_i = (\beta_{i,1}, \beta_{i,2}, \cdots, \beta_{i,n}) \in \mathbb{N}^n$ and $h_{\boldsymbol{\beta}_i} \in \mathbb{R}$ is corresponding real-valued coefficient.

*Illustration* We consider an model updating example where a scalar stiffness updating variable $\theta$ represents the relative change from the initial/nominal value of a stiffness parameter. Another variable for the optimization problem is the unmeasured 4th entry in the first mode shape vector, $\boldsymbol{\psi}_{1,u} = \psi_{1,4}$ (abbreviated as $\psi_4$ herein). They constitute the optimization vector variable $\mathbf{x} = (x_1, x_2)^\mathrm{T} = (\theta, \psi_4)^\mathrm{T}$, i.e. $n = 2$. After plugging in the numerical values of the example structure, the corresponding optimization problem from the modal dynamic residual approach described in Eq. (3) is found as follows. Because the objective function in Eq. (3) is the square of a $\mathcal{L}_2$-norm, the highest-degree term in the expanded polynomial is $200\theta^2 \psi_4^2$, i.e. degree $d = 4$. For generality and for illustration, in the objective function we also list the six monomials with coefficient 0 (i.e. non-existing) that have degree less than or equal to $d$.

$$\begin{aligned} \underset{\theta,\psi_4}{\text{minimize}} \quad & f(\mathbf{x}) = 229.584 + 427.670\theta - 403.687\psi_4 + 200\theta^2 - 803.687\theta\psi_4 + \\ & \quad 177.455\psi_4^2 + 0 \cdot \theta^3 - 400\theta^2\psi_4 + 376.017\theta\psi_4^2 + 0 \cdot \psi_4^3 + 0 \cdot \theta^4 + \\ & \quad 0 \cdot \theta^3\psi_4 + 200\theta^2\psi_4^2 + 0 \cdot \theta\psi_4^3 + 0 \cdot \psi_4^4 \\ \text{subject to} \quad & g_1(\mathbf{x}) = (1-\theta)(1+\theta) \geq 0 \\ & g_2(\mathbf{x}) = (2-\psi_4)(2+\psi_4) \geq 0 \end{aligned} \tag{6}$$

Note that to apply the SOS method, the box constraints from Eq. (3) are equivalently rewritten into polynomial functions $g_1(\mathbf{x})$ and $g_2(\mathbf{x})$. With $L_\theta = -1$ (i.e. $-100\%$) and $U_\theta = 1$, $L_\theta \leq \theta \leq U_\theta$ is equivalent to $g_1(\mathbf{x}) \geq 0$. With $\boldsymbol{\psi}_{1,m}$ normalized to have maximum magnitude of 1, the coefficients in $g_2(\mathbf{x})$ bound the unmeasured entry as $\boldsymbol{\psi}_{1,u} = \psi_{1,4} \in [-2, 2]$. If relaxation to the bounds is needed, we can easily change the constant coefficients in $g_2(\mathbf{x})$. In this optimization problem, the degrees of the objective function and inequality constraints are $d = 4$ from the term $200\theta^2\psi_4^2$ in $f(\mathbf{x})$, $e_1 = 2$ from the term $-\theta^2$ in $g_1(\mathbf{x})$, and $e_2 = 2$ from the term $-\psi_4^2$ in $g_2(\mathbf{x})$, respectively. The total number of different power vectors $\boldsymbol{\alpha} \in \mathbb{N}^n$ with the vector sum $\sum_{i=1}^n \alpha_i \leq d$, i.e. total number of base monomials in the objective function in Eq. (6), equals $\binom{n+d}{n} = \binom{2+4}{2} = 15$. Taking the monomial $376.017\theta\psi_4^2 = 376.017 x_1 x_2^2$ for example, the power vector $\boldsymbol{\alpha} = (1,2)$ and the coefficient $c_{\boldsymbol{\alpha}} = 376.017$. ∎

Define the feasible set of the optimization problem in Eq. (5) as $\boldsymbol{\Omega} = \{\mathbf{x} \in \mathbb{R}^n | g_i(\mathbf{x}) \geq 0, i = 1, 2, \cdots, l\}$. If $f^* = f(\mathbf{x}^*)$ is the global minimum value of the problem, $f(\mathbf{x}) - f^*$ is nonnegative for all $\mathbf{x} \in \boldsymbol{\Omega}$. Therefore, the original optimization problem in Eq. (5) can be equivalently reformulated as finding the maximum lower bound of a scalar objective $\gamma$:

$$\begin{aligned} \underset{\gamma}{\text{maximize}} \quad & \gamma \\ \text{subject to} \quad & f(\mathbf{x}) - \gamma \geq 0, \quad \forall \mathbf{x} \in \boldsymbol{\Omega} \end{aligned} \tag{7}$$

The optimal objective value, $\gamma^*$, is intended to approach $f^*$. Despite the general nonconvexity of $\Omega$ on $\mathbf{x}$, for each (fixed) $\mathbf{x} \in \Omega$, $f(\mathbf{x}) - \gamma$ is an affine function of $\gamma$, and thus $f(\mathbf{x}) - \gamma \geq 0$ is a convex constraint on $\gamma$. The feasible set of $\gamma$ in Eq. (7) is therefore the intersection of infinite number of convex constraints on $\gamma$. As a result, the set is still convex on $\gamma$, and the optimization problem is convex [26].

As it remains a hard problem to test that a polynomial $f(\mathbf{x}) - \gamma$ is nonnegative for all $\mathbf{x} \in \Omega$, we attempt to decompose $f(\mathbf{x}) - \gamma$ as sum of squares (SOS) over $\Omega$, which immediately implies nonnegativity (as shown in Eq. (1)). According to Lasserre [23], a sufficient and more solvable condition for the nonnegativity of $f(\mathbf{x}) - \gamma$ over $\Omega$ is that there exist SOS polynomials $s_i(\mathbf{x}) = \mathbf{z}_i^\mathrm{T} \mathbf{Q}_i \mathbf{z}_i$, and $\mathbf{Q}_i \succcurlyeq 0, i = 0, 1, \cdots, l$, that satisfy the following SOS decomposition of $f(\mathbf{x}) - \gamma$:

$$f(\mathbf{x}) - \gamma = s_0(\mathbf{x}) + \sum_{i=1}^{l} s_i(\mathbf{x}) g_i(\mathbf{x}) \tag{8}$$

Recall that polynomials $f(\mathbf{x})$ and $g_i(\mathbf{x})$ have degree $d$ and $e_i, i = 1, 2, \cdots, l$, respectively. In order to make the equality in Eq. (8) hold, we express both sides as a polynomial with degree of $2t$, where $t$ is the smallest integer such that $2t \geq \max(d, e_1, \cdots, e_l)$. In other words, all polynomials $f(\mathbf{x})$ and $g_i(\mathbf{x}), i = 1, 2, \cdots, l$ have degree no more than $2t$. When expressing $f(\mathbf{x}) - \gamma$ at the left hand side as a polynomial with degree of $2t$, if it turns out that $d < 2t$, the monomials with degree larger than $d$ are simply assigned with zero coefficient.

On the right-hand side of Eq. (8), to ensure first the degree of the polynomial $s_0(\mathbf{x}) = \mathbf{z}_0^\mathrm{T} \mathbf{Q}_0 \mathbf{z}_0$ is no more than $2t$, we define $\mathbf{z}_0$ to represent all the base monomials of degree $t \in \mathbb{N}$ or lower:

$$\mathbf{z}_0 = (1, x_1, x_2, \cdots, x_n, x_1^2, x_1 x_2, \cdots x_{n-1} x_n^{t-1}, x_n^t)^\mathrm{T} \tag{9}$$

The length of $\mathbf{z}_0$ is $n_{\mathbf{z}_0} = \binom{n+t}{n}$, and $\mathbf{Q}_0 \in \mathbb{S}_+^{n_{\mathbf{z}_0}}$, i.e. the set of symmetric positive semi-definite $n_{\mathbf{z}_0} \times n_{\mathbf{z}_0}$ matrices. Likewise, to ensure each product $s_i(\mathbf{x}) g_i(\mathbf{x}), i = 1, 2, \cdots, l$, has degree no more than $2t$, $\mathbf{z}_i$ is defined as the vector including all the monomials of degree $t - \tilde{e}_i$ or lower, where $\tilde{e}_i = \lceil e_i/2 \rceil$ is the smallest integer such that $\tilde{e}_i \geq e_i/2$. As a result, the length of $\mathbf{z}_i$ is $n_{\mathbf{z}_i} = \binom{n+t-\tilde{e}_i}{n}$ and $\mathbf{Q}_i \in \mathbb{S}_+^{n_{\mathbf{z}_i}}$. Because the degree of $g_i(\mathbf{x})$ is $e_i$, and the degree of $s_i(\mathbf{x})$ is $2(t - \tilde{e}_i)$, the product $s_i(\mathbf{x}) g_i(\mathbf{x})$ has degree of $2t - 2\tilde{e}_i + e_i$, which is guaranteed to be no more than $2t$. In summary, the SOS decomposition in Eq. (8) guarantees $f(\mathbf{x}) - \gamma$ to be nonnegative for all $\mathbf{x} \in \Omega = \{\mathbf{x} \in \mathbb{R}^n | g_i(\mathbf{x}) \geq 0, i = 1, 2, \cdots, l\}$, which is equivalent to the constraint in Eq. (7).

*Illustration – continued* We start with finding the degrees $2t = 4 = \max(d, e_1, e_2)$, $\tilde{e}_1 = \lceil e_1/2 \rceil = 1, \tilde{e}_2 = \lceil e_2/2 \rceil = 1$. The vector lengths are determined as $n_{\mathbf{z}_0} = \binom{2+2}{2} = 6$ and $n_{\mathbf{z}_1} = n_{\mathbf{z}_2} = \binom{2+2-1}{2} = 3$. The base monomial vectors are then defined as:

$$\mathbf{z}_0 = (1, \theta, \psi_4, \theta \psi_4, \theta^2, \psi_4^2)^\mathrm{T}, \mathbf{z}_1 = (1, \theta, \psi_4)^\mathrm{T}, \mathbf{z}_2 = (1, \theta, \psi_4)^\mathrm{T} \tag{10}$$

The nonnegativity condition of $f(\mathbf{x}) - \gamma$ over $\Omega$ is that there exist positive semidefinite matrices $\mathbf{Q}_0, \mathbf{Q}_1, \mathbf{Q}_2$ such that $f(\mathbf{x}) - \gamma = \mathbf{z}_0^\mathrm{T} \mathbf{Q}_0 \mathbf{z}_0 + \mathbf{z}_1^\mathrm{T} \mathbf{Q}_1 \mathbf{z}_1 g_1(\mathbf{x}) + \mathbf{z}_2^\mathrm{T} \mathbf{Q}_2 \mathbf{z}_2 g_2(\mathbf{x})$. ∎

Using this SOS sufficient condition of polynomial nonnegativity, with $s_i(\mathbf{x}) = \mathbf{z}_i^T \mathbf{Q}_i \mathbf{z}_i$ and $\mathbf{Q}_i \in \mathbb{S}_+^{n_{z_i}}$, $i = 0, 1, 2, \cdots, l$, the optimization problem described in Eq. (7) can be relaxed to a semi-definite programming (SDP) problem:

$$\begin{aligned}
\underset{\gamma, \mathbf{Q}_i}{\text{maximize}} \quad & \gamma \\
\text{subject to} \quad & f(\mathbf{x}) - \gamma = \mathbf{z}_0^T \mathbf{Q}_0 \mathbf{z}_0 + \sum_{i=1}^{l} (\mathbf{z}_i^T \mathbf{Q}_i \mathbf{z}_i) g_i(\mathbf{x}) \\
& \mathbf{Q}_i \succcurlyeq 0, \ i = 0, 1, \cdots, l
\end{aligned} \quad (11)$$

Note that the identity in Eq. (11) is an equality constraint that holds for arbitrary $\mathbf{x}$, which essentially says two sides of the equation should have the same coefficient $c_\alpha$ for the same base monomial $\mathbf{x}^\alpha$. At the left-hand side, variable $\gamma$ is contained in the constant coefficient. This equality constraint is effectively a group of affine equality constraints on the entries of $\mathbf{Q}_i$ and $\gamma$. The total number of such affine equality constraints is $\binom{n+2t}{n}$, which equals to the number of all base monomials of $\mathbf{x} \in \mathbb{R}^n$ with degree less than or equal to $2t$ (and equals the total number of different power vectors $\boldsymbol{\alpha} \in \mathbb{N}^n$ with the vector sum $\sum_{i=1}^{n} \alpha_i \leq 2t$).

*Illustration – continued*      To apply SOS optimization method, we introduce variables $\gamma$, $\mathbf{Q}_0 = [Q_{ij}^{(0)}]_{6\times6}$, $\mathbf{Q}_1 = [Q_{ij}^{(1)}]_{3\times3}$, and $\mathbf{Q}_2 = [Q_{ij}^{(2)}]_{3\times3}$. In this illustration, because $d > e_1$ and $d > e_2$, we have $2t = d$. The identity on $\mathbf{x}$ generates $\binom{n+2t}{n} = \binom{2+4}{2} = 15$ equality constraints, which correspond to the 15 monomials. All monomial coefficients of $f(\theta, \psi_4) - \gamma$ can be directly read from Eq. (6), except that the constant term becomes $229.584 - \gamma$. The problem in Eq. (6) is then relaxed to the following SDP problem.

$$\begin{aligned}
&\underset{\gamma, \mathbf{Q}_0, \mathbf{Q}_1, \mathbf{Q}_2}{\text{maximize}} \quad \gamma \\
&\text{subject to} \quad \mathbf{Q}_0 = [Q_{ij}^{(0)}]_{6\times6} \succcurlyeq 0, \quad\quad \mathbf{Q}_1 = [Q_{ij}^{(1)}]_{3\times3} \succcurlyeq 0, \quad\quad \mathbf{Q}_2 = [Q_{ij}^{(2)}]_{3\times3} \succcurlyeq 0, \\
&\text{Constant:} \quad 229.584 - \gamma = Q_{11}^{(0)} + Q_{11}^{(1)} + 4Q_{11}^{(2)}, \quad\quad 427.670\theta: \ 427.670 = 2Q_{12}^{(0)} + 2Q_{12}^{(1)} + 8Q_{12}^{(2)}, \\
&-403.687\psi_4: \ -403.687 = 2Q_{13}^{(0)} + 2Q_{13}^{(1)} + 8Q_{13}^{(2)}, \quad\quad 200\theta^2: \ 200 = Q_{22}^{(0)} + 2Q_{15}^{(0)} - Q_{11}^{(1)} + Q_{22}^{(1)} + 4Q_{22}^{(2)}, \\
&-803.687\theta\psi_4: \ -803.687 = 2Q_{14}^{(0)} + 2Q_{23}^{(0)} + 2Q_{23}^{(1)} + 8Q_{23}^{(2)}, \quad 177.455\psi_4^2: \ 177.455 = Q_{33}^{(0)} + 2Q_{16}^{(0)} + Q_{33}^{(1)} - Q_{11}^{(2)} + 4Q_{33}^{(2)}, \\
&0\cdot\theta^3: \ 0 = 2Q_{25}^{(0)} - 2Q_{12}^{(1)}, \quad\quad -400\theta^2\psi_4: \ -400 = 2Q_{24}^{(0)} + 2Q_{35}^{(0)} - 2Q_{13}^{(1)}, \\
&376.017\theta\psi_4^2: \ 376.017 = 2Q_{26}^{(0)} + 2Q_{34}^{(0)} - 2Q_{12}^{(2)}, \quad\quad 0\cdot\psi_4^3: \ 0 = 2Q_{36}^{(0)} - 2Q_{13}^{(2)}, \\
&0\cdot\theta^4: \ 0 = Q_{55}^{(0)} - Q_{22}^{(1)}, \quad\quad 0\cdot\theta^3\psi_4: \ 0 = 2Q_{45}^{(0)} - 2Q_{23}^{(1)}, \\
&200\theta^2\psi_4^2: \ 200 = Q_{44}^{(0)} + 2Q_{56}^{(0)} - Q_{33}^{(1)} - Q_{22}^{(2)}, \quad\quad 0\cdot\theta\psi_4^3: \ 0 = 2Q_{46}^{(0)} - 2Q_{23}^{(2)}, \\
&0\cdot\psi_4^4: \ 0 = Q_{66}^{(0)} - Q_{33}^{(2)},
\end{aligned} \quad (12)$$

The equality constraints in Eq. (12) are obtained by expanding both sides of the equality constraints in Eq. (11) with arbitrary $\mathbf{x}$ into sum of monomials, and then setting the corresponding monomial coefficients to be identical. For every equality constraint, the corresponding monomial is listed to the left of it. Take the first equality constraint, $229.584 - \gamma = Q_{11}^{(0)} + Q_{11}^{(1)} + 4Q_{11}^{(2)}$, for example. On the left-hand side, $229.584 - \gamma$ is the constant term of the polynomial $f(\mathbf{x}) - \gamma$. On the right-hand side, $Q_{11}^{(0)} + Q_{11}^{(1)} + 4Q_{11}^{(2)}$ is the constant term upon expanding the polynomial $\mathbf{z}_0^T \mathbf{Q}_0 \mathbf{z}_0 + \sum_{i=1}^{l} (\mathbf{z}_i^T \mathbf{Q}_i \mathbf{z}_i) g_i(\mathbf{x})$ into sum of monomials. The coefficient 4 in front of $Q_{11}^{(2)}$ comes from the constant term in $g_2(\mathbf{x})$ from Eq. (6). ∎

In order to conveniently derive the dual problem later, the group of equality constraints in Eq. (11) are now equivalently and explicitly rewritten using constant selection matrices $\mathbf{A}_\alpha$ and $\mathbf{B}_{i,\alpha}$ ($i = 1, 2, \cdots, l$). Let $\mathbb{S}$ denote the set of the real-valued symmetric matrices. We assign matrix $\mathbf{A}_\alpha \in \mathbb{S}^{n_{z_0}}$ to have value 1 for entries where $\mathbf{x}^\alpha$ appears in the square matrix $\mathbf{z}_0 \mathbf{z}_0^T \in \mathbb{S}^{n_{z_0}}$ and value 0 otherwise (recall $\mathbf{z}_0$ from Eq. (9)). In other words, $\mathbf{A}_\alpha$ selects entries with $\mathbf{x}^\alpha$ from the matrix $\mathbf{z}_0 \mathbf{z}_0^T$. Similarly, for each $g_i(\mathbf{x}) = \sum_{\boldsymbol{\beta}_i} h_{\boldsymbol{\beta}_i} \mathbf{x}^{\boldsymbol{\beta}_i}$, matrix $\mathbf{B}_{i,\alpha} \in \mathbb{S}^{n_{z_i}}$ has value $h_{\boldsymbol{\beta}_i}$ in entries where $\mathbf{x}^\alpha$ appears in the square matrix $\mathbf{x}^{\boldsymbol{\beta}_i} \mathbf{z}_i \mathbf{z}_i^T$ and value 0 otherwise. Using operator $\langle \cdot, \cdot \rangle$ to denote the matrix inner product, the optimization problem in Eq. (11) can be equivalently rewritten as follows, with $\gamma$ and matrices $\mathbf{Q}_i \in \mathbb{S}_+^{n_{z_i}}$, $i = 0, 1, 2, \cdots, l$, as optimization variables.

$$\begin{aligned}
\underset{\gamma, \mathbf{Q}_i}{\text{maximize}} \quad & \gamma \\
\text{subject to} \quad & \langle \mathbf{A}_0, \mathbf{Q}_0 \rangle + \sum_{i=1}^{l} \langle \mathbf{B}_{i,0}, \mathbf{Q}_i \rangle = c_0 - \gamma \\
& \langle \mathbf{A}_\alpha, \mathbf{Q}_0 \rangle + \sum_{i=1}^{l} \langle \mathbf{B}_{i,\alpha}, \mathbf{Q}_i \rangle = c_\alpha \qquad \text{for all } \boldsymbol{\alpha} \neq \mathbf{0} \\
& \mathbf{Q}_i \succcurlyeq 0 \qquad\qquad\qquad\qquad\qquad i = 0, 1, 2, \cdots, l
\end{aligned} \qquad (13)$$

Here $\gamma$ and $\mathbf{Q}_i$ are optimization variables; $c_\alpha$ is the real coefficient for the monomial $\mathbf{x}^\alpha$ in $f(\mathbf{x})$; $\mathbf{A}_\alpha$ and $\mathbf{B}_{i,\alpha}$ are constant matrices. The total number of equality constraints is still $\binom{n+2t}{n}$. This equals the number of base monomials of $\mathbf{x} \in \mathbb{R}^n$ with degree of $2t$ or lower, i.e. the number of different power vectors $\boldsymbol{\alpha} \in \mathbb{N}^n$ (including the zero vector $\mathbf{0}$) with the vector sum $\sum_{i=1}^{n} \alpha_i \leq 2t$. The problem in Eq. (13) can be easily input into standard convex optimization solvers.

*Illustration – continued*   Taking the monomial $\theta^2$ in Eq. (6) for example, the power vector $\boldsymbol{\alpha} = (2,0)$, the coefficient $c_\alpha = 200$, $\mathbf{x}^\alpha = x_1^2 x_2^0 = \theta^2$. Recalling that $\mathbf{z}_0 = (1, \theta, \psi_4, \theta\psi_4, \theta^2, \psi_4^2)^T$, $\theta^2$ appears at entries $(1,5)$, $(2,2)$, and $(5,1)$ in matrix $\mathbf{z}_0 \mathbf{z}_0^T$. Therefore, $\mathbf{A}_\alpha$ has value 1 at these three entries and value 0 at all other entries. We then consider $\mathbf{z}_1 = (1, \theta, \psi_4)^T$ and $g_1(\mathbf{x}) = 1 - \theta^2$. The two power vectors in $g_1(\mathbf{x})$ are $\boldsymbol{\beta}_{1,1} = (0,0)$ corresponding to the constant term $\mathbf{x}^{\boldsymbol{\beta}_{1,1}} = 1$, and $\boldsymbol{\beta}_{1,2} = (2,0)$ corresponding to $\mathbf{x}^{\boldsymbol{\beta}_{1,2}} = \theta^2$. Coefficients $h_{\boldsymbol{\beta}_{1,1}} = 1$ and $h_{\boldsymbol{\beta}_{1,2}} = -1$. Because $\mathbf{x}^\alpha = \theta^2$ appears at entry $(2,2)$ in matrix $\mathbf{x}^{\boldsymbol{\beta}_{1,1}} \cdot \mathbf{z}_1 \mathbf{z}_1^T = 1 \cdot \mathbf{z}_1 \mathbf{z}_1^T$ and entry $(1,1)$ in matrix $\mathbf{x}^{\boldsymbol{\beta}_{1,2}} \cdot \mathbf{z}_1 \mathbf{z}_1^T = \theta^2 \cdot \mathbf{z}_1 \mathbf{z}_1^T$, $\mathbf{B}_{1,\alpha}$ has value $h_{\boldsymbol{\beta}_{1,1}} = 1$ at entry $(2,2)$, value $h_{\boldsymbol{\beta}_{1,2}} = -1$ at entry $(1,1)$, and value 0 at all other entries. Similarly, as $\mathbf{z}_2 = (1, \theta, \psi_4)^T$ and $g_2(\mathbf{x}) = 4 - \psi_4^2$, we have $\boldsymbol{\beta}_{2,1} = (0,0)$ and $\boldsymbol{\beta}_{2,2} = (0,2)$ and coefficients $h_{\boldsymbol{\beta}_{2,1}} = 4$ and $h_{\boldsymbol{\beta}_{2,2}} = -1$. Because $\mathbf{x}^\alpha = \theta^2$ appears at entry $(2,2)$ in matrix $\mathbf{x}^{\boldsymbol{\beta}_{2,1}} \cdot \mathbf{z}_2 \mathbf{z}_2^T = 1 \cdot \mathbf{z}_2 \mathbf{z}_2^T$ and nowhere in matrix $\mathbf{x}^{\boldsymbol{\beta}_{2,2}} \cdot \mathbf{z}_2 \mathbf{z}_2^T = \psi_4^2 \cdot \mathbf{z}_2 \mathbf{z}_2^T$, $\mathbf{B}_{2,\alpha}$ has value $h_{\boldsymbol{\beta}_{2,1}} = 4$ at entry $(2,2)$ and value 0 at all other entries. For $\boldsymbol{\alpha} = (2,0)$, the selection matrices $\mathbf{A}_\alpha$, $\mathbf{B}_{1,\alpha}$, and $\mathbf{B}_{2,\alpha}$ are shown as:

$$\mathbf{A}_\alpha = \begin{bmatrix} 0 & 0 & 0 & 0 & 1 & 0 \\ 0 & 1 & 0 & 0 & 0 & 0 \\ 0 & 0 & 0 & 0 & 0 & 0 \\ 0 & 0 & 0 & 0 & 0 & 0 \\ 1 & 0 & 0 & 0 & 0 & 0 \\ 0 & 0 & 0 & 0 & 0 & 0 \end{bmatrix} \qquad \mathbf{B}_{1,\alpha} = \begin{bmatrix} -1 & 0 & 0 \\ 0 & 1 & 0 \\ 0 & 0 & 0 \end{bmatrix} \qquad \mathbf{B}_{2,\alpha} = \begin{bmatrix} 0 & 0 & 0 \\ 0 & 4 & 0 \\ 0 & 0 & 0 \end{bmatrix}$$

Noticing that matrices $\mathbf{Q}_0$, $\mathbf{Q}_1$, and $\mathbf{Q}_2$ are symmetric, the equality on coefficient of $\theta^2$ can be written as $Q_{22}^{(0)} + 2Q_{15}^{(0)} - Q_{11}^{(1)} + Q_{22}^{(1)} + 4Q_{22}^{(2)} = 200$, which is the same as the expression in Eq. (12). Other equality constraints can be formulated in a similar way. ∎

Upon solving the optimization problem in Eq. (13) through SOS relaxation, the best lower bound, i.e. the largest $\gamma^*$ such that $\gamma^* \leq f^*$, of the objective function in Eq. (5) is obtained. For most practical applications, the lower bound obtained by SOS relaxation usually coincides with the optimal value of the objective function, i.e. $\gamma^* = f^*$ [21].

To summarize the optimization procedure, FIGURE 1 shows the flow chat of the procedure. First, the problem of minimizing a polynomial $f(\mathbf{x})$ over a set $\mathbf{\Omega}$ (Eq. (5)) is equivalently reformulated as finding the best lower bound $\gamma^*$ of $f(\mathbf{x})$ over the set $\mathbf{\Omega}$ (Eq. (7)). Second, the condition that $f(\mathbf{x}) - \gamma \geq 0$ over set $\mathbf{\Omega}$ is relaxed to a more easily solvable condition that $f(\mathbf{x}) - \gamma$ has an SOS decomposition over set $\mathbf{\Omega}$ (Eq. (11) and Eq. (13)).

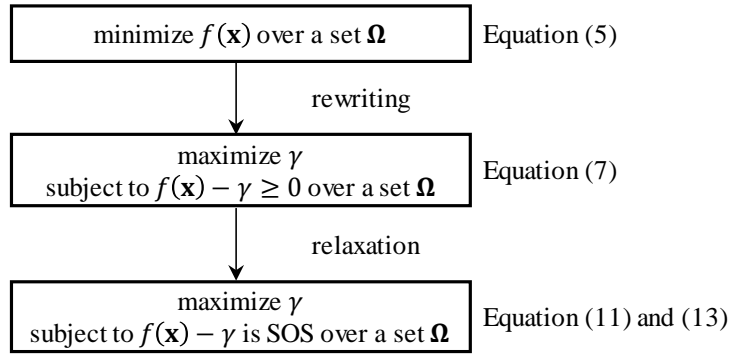

FIGURE 1 Flow chat of the optimization procedure

## 3.2 Dual problem

To accomplish model updating, only finding the lower bound or the optimal value of the objective function ($f^*$) is not enough. The minimizer of the objective function, $\mathbf{x}^*$ such that $f(\mathbf{x}^*) = f^*$, needs to be computed because $\mathbf{x}^*$ contains optimal values of the stiffness updating variables. Fortunately, the minimizer can be easily extracted from the solution of the dual problem of Eq. (13) [24, 27]. Define the dual variables, a.k.a. Lagrangian multiplier vectors $\mathbf{y} \in \mathbb{R}^{\binom{n+2t}{n}}$ and matrices $\mathbf{U}_i \in \mathbb{S}^{n_{z_i}}$, $i = 0, 1, 2, \cdots, l$. Dual variable $y_i$, $i = 1, 2, \cdots, \binom{n+2t}{n}$, is associated with the $i$-th equality constraint in Eq. (13). Variable $\mathbf{U}_i$, $i = 0, 1, 2, \cdots, l$, is associated with the $i$-th inequality constraint in Eq. (13). The Lagrangian for the primal problem in Eq. (13) is:

$$\mathcal{L}(\gamma, \mathbf{Q}_i, \mathbf{y}, \mathbf{U}_i) = \gamma + \left(c_0 - \gamma - \langle \mathbf{A}_0, \mathbf{Q}_0 \rangle - \sum_{i=1}^{l} \langle \mathbf{B}_{i,0}, \mathbf{Q}_i \rangle\right) y_0 + \sum_{\alpha \neq 0} \left(c_\alpha - \langle \mathbf{A}_\alpha, \mathbf{Q}_0 \rangle - \sum_{i=1}^{l} \langle \mathbf{B}_{i,\alpha}, \mathbf{Q}_i \rangle\right) y_\alpha + \sum_{i=0}^{l} \langle \mathbf{U}_i, \mathbf{Q}_i \rangle$$
$$= \sum_\alpha c_\alpha y_\alpha + \gamma(1 - y_0) + \langle \mathbf{U}_0 - \sum_\alpha \mathbf{A}_\alpha y_\alpha, \mathbf{Q}_0 \rangle + \sum_{i=1}^{l} \langle \mathbf{U}_i - \sum_\alpha \mathbf{B}_{i,\alpha} y_\alpha, \mathbf{Q}_i \rangle \quad (14)$$

The dual function is then formed as the supremum of the Lagrangian with respect to primal variables $\gamma$ and $\mathbf{Q}_i$. Because $\mathcal{L}$ is affine on $\gamma$ and $\mathbf{Q}_i$, the dual function is found as:

$$\mathcal{D}(\mathbf{y}, \mathbf{U}_i) = \sup_{\gamma, \mathbf{Q}_i} \mathcal{L}(\gamma, \mathbf{Q}_i, \mathbf{y}, \mathbf{U}_i) = \begin{cases} \sum_\alpha c_\alpha y_\alpha & \text{if } y_0 = 1, \mathbf{U}_0 = \sum_\alpha \mathbf{A}_\alpha y_\alpha, \mathbf{U}_i = \sum_\alpha \mathbf{B}_{i,\alpha} y_\alpha \\ +\infty & \text{otherwise} \end{cases} \quad (15)$$

As a result, the dual problem of the primal problem in Eq. (13) is written as:

$$\begin{aligned} \underset{\mathbf{y}}{\text{minimize}} \quad & \sum_\alpha c_\alpha y_\alpha \\ \text{subject to} \quad & y_0 = 1 \\ & \mathbf{U}_0 = \sum_\alpha \mathbf{A}_\alpha y_\alpha \succcurlyeq 0 \\ & \mathbf{U}_i = \sum_\alpha \mathbf{B}_{i,\alpha} y_\alpha \succcurlyeq 0 \quad i = 1, 2, \cdots, l \end{aligned} \quad (16)$$

where $\mathbf{y} = \{y_\alpha\} \in \mathbb{R}^{\binom{n+2t}{n}}$ is the optimization vector variable. $\mathbf{A}_\alpha$ and $\mathbf{B}_{i,\alpha}$ ($i = 1, 2, \cdots, l$) are matrices described in Section 3.1.

It has been shown that if the optimal value of the original problem (Eq. (5)) and the SOS primal problem (Eq. (13)) coincide with each other, the optimal solution of the SOS dual problem can be calculated as:

$$\mathbf{y}^* = (1, x_1^*, \cdots, x_n^*, (x_1^*)^2, x_1^* x_2^* \cdots, (x_n^*)^{2t})^{\text{T}} \quad (17)$$

where the entries correspond to the monomials $\mathbf{x}^\alpha$ [23]. Finally, the optimal solution $\mathbf{x}^*$ for the original problem in Eq. (5) is now easily extracted from $\mathbf{y}^*$, as the 2nd through the $(n+1)$-th entries. We refer interested readers to Henrion and Lasserre [25] for details of the minimizer extracting technique. Since practical SDP solvers, such as SeDuMi [28], simultaneously solve both primal and dual problems, the optimal solution $\mathbf{x}^*$ can be computed efficiently.

*Illustration – continued*    For the dual problem, the optimization vector variable is defined as follows. Each subscript $(\cdot,\cdot)$ denotes a power vector $\boldsymbol{\alpha} \in \mathbb{N}^2$ such that $\sum_{i=1}^{2} \alpha_i \leq 2t = 4$.

$$\mathbf{y} = (y_{(0,0)}, y_{(1,0)}, y_{(0,1)}, y_{(2,0)}, y_{(1,1)}, y_{(0,2)}, y_{(3,0)}, y_{(2,1)}, y_{(1,2)}, y_{(0,3)}, y_{(4,0)}, y_{(3,1)}, y_{(2,2)}, y_{(1,3)}, y_{(0,4)})^{\text{T}}$$

And the dual problem is formulated as:

$$\begin{aligned} \underset{\mathbf{y}}{\text{minimize}} \quad & 229.584 y_{(0,0)} + 427.670 y_{(1,0)} - 403.687 y_{(0,1)} + 200 y_{(2,0)} - 803.687 y_{(1,1)} + \\ & 177.455 y_{(0,2)} + 0 y_{(3,0)} - 400 y_{(2,1)} + 376.017 y_{(1,2)} + 0 y_{(0,3)} + 0 y_{(4,0)} + \\ & 0 y_{(3,1)} + 200.00 y_{(2,2)} + 0 y_{(1,3)} + 0 y_{(0,4)} \\ \text{subject to} \quad & y_{(0,0)} = 1 \\ & \mathbf{U}_0 = \begin{bmatrix} y_{(0,0)} & y_{(1,0)} & y_{(0,1)} & y_{(1,1)} & y_{(2,0)} & y_{(0,2)} \\ y_{(1,0)} & y_{(2,0)} & y_{(1,1)} & y_{(2,1)} & y_{(3,0)} & y_{(1,2)} \\ y_{(0,1)} & y_{(1,1)} & y_{(0,2)} & y_{(1,2)} & y_{(2,1)} & y_{(0,3)} \\ y_{(1,1)} & y_{(2,1)} & y_{(1,2)} & y_{(2,2)} & y_{(3,1)} & y_{(1,3)} \\ y_{(2,0)} & y_{(3,0)} & y_{(2,1)} & y_{(3,1)} & y_{(4,0)} & y_{(2,2)} \\ y_{(0,2)} & y_{(1,2)} & y_{(0,3)} & y_{(1,3)} & y_{(2,2)} & y_{(0,4)} \end{bmatrix} \succcurlyeq 0 \\ & \mathbf{U}_1 = \begin{bmatrix} y_{(0,0)} - y_{(2,0)} & y_{(1,0)} - y_{(3,0)} & y_{(0,1)} - y_{(2,1)} \\ y_{(1,0)} - y_{(3,0)} & y_{(2,0)} - y_{(4,0)} & y_{(1,1)} - y_{(3,1)} \\ y_{(0,1)} - y_{(2,1)} & y_{(1,1)} - y_{(3,1)} & y_{(0,2)} - y_{(2,2)} \end{bmatrix} \succcurlyeq 0 \end{aligned} \quad (18)$$

$$\mathbf{U}_2 = \begin{bmatrix} 4y_{(0,0)} - y_{(0,2)} & 4y_{(1,0)} - y_{(1,2)} & 4y_{(0,1)} - y_{(0,3)} \\ 4y_{(1,0)} - y_{(1,2)} & 4y_{(2,0)} - y_{(2,2)} & 4y_{(1,1)} - y_{(1,3)} \\ 4y_{(0,1)} - y_{(0,3)} & 4y_{(1,1)} - y_{(1,3)} & 4y_{(0,2)} - y_{(0,4)} \end{bmatrix} \succcurlyeq 0$$

Recall the $\mathbf{A}_\alpha$, $\mathbf{B}_{1,\alpha}$, and $\mathbf{B}_{2,\alpha}$ matrices provided for $\alpha = (2, 0)$ in the last illustration. Consequently, dual variable $y_{(2,0)}$ appears in the same entries of $\mathbf{U}_0$ and $\mathbf{A}_\alpha$, the same entries of $\mathbf{U}_1$ and $\mathbf{B}_{1,\alpha}$, and the same entries of $\mathbf{U}_2$ and $\mathbf{B}_{2,\alpha}$.  ∎

As all the functions in modal dynamic residual approach in Eq. (3) are polynomials, the primal and dual problems of SOS optimization can be directly implemented. In this way, the modal dynamic residual approach is recast as a convex problem. CVX, an MATLAB interface package for defining convex problems is used to solve the problems in this paper [29]. SeDuMi is adopted as the underlying solver.

## 4  Validation examples

### 4.1  Numerical simulation

To validate the proposed sum of squares (SOS) method for model updating, a four-story shear frame is first simulated (FIGURE 2). In the initial model, the nominal weight and inter-story stiffness values of all floors are set as 12.060 lb and 10 lbf/in, respectively. To construct an "as-built" structure, the stiffness value of the fourth story is reduced by 10% to 9 lbf/in, as shown in FIGURE 2. Modal properties of the "as-built" structure are directly used as "experimental" modal properties. It is assumed that only the first three floors are instrumented with sensors and only the first mode is "measured" and available for model updating. The first resonance frequency and the "measured" (instrumented) three mode shape entries are $\omega_1 = 6.196 \text{ rad/s}$ and $\boldsymbol{\psi}_\text{m} = (0.395, 0.742, 1.000)^\text{T}$, respectively. Recalling notations in Eq. (3), here $n_\text{modes} = 1$. The vector of measured three entries of the first mode shape $\boldsymbol{\psi}_{1,\text{m}}$ is abbreviated as $\boldsymbol{\psi}_\text{m}$, while the maximum among the three entries are normalized to be 1.

To make 3D graphical illustration possible for the nonconvex objective function, the optimization variables include the stiffness parameter change only of the fourth floor, denoted as $\theta$, and the fourth entry in the mode shape vector, denoted as $\psi_{1,4}$ and abbreviated as $\psi_4$. Here the variable $\theta$ represents the relative change of $k_4$ from the initial nominal value of 10 lbf/in, i.e. $\theta = (k_4 - 10)/10$. In other words, it is assumed $k_1$, $k_2$, and $k_3$ do not require updating; we know the change happens with $k_4$ but need to identify/update how much the change is. The value of $\psi_4$ is obviously influenced by the previous normalization in $\boldsymbol{\psi}_\text{m}$. With only two optimization variables, $\theta$ and $\psi_4$, the nonconvex objective function can be written as:

$$\begin{aligned}
\underset{\theta, \psi_4}{\text{minimize}} \quad & f(\mathbf{x}) = f(\theta, \psi_4) := \left\| \left[ (\mathbf{K}_0 + \theta \mathbf{K}_{0,4}) - \omega_1^2 \mathbf{M} \right] \begin{Bmatrix} \boldsymbol{\Psi}_\text{m} \\ \psi_4 \end{Bmatrix} \right\|_2^2 \\
\text{subject to} \quad & g_1(\mathbf{x}) = (1-\theta)(1+\theta) \geq 0 \\
& g_2(\mathbf{x}) = (2-\psi_4)(2+\psi_4) \geq 0
\end{aligned} \quad (19)$$

Using the adopted numeric values and expanding $f(\mathbf{x})$ into polynomial form, the model updating problem is formulated as shown in the illustration in Eq. (6). This formulation equivalently sets the lower bounds of optimization variables as $L_\theta = -1$ (meaning lowest possible stiffness value is 0) and $L_{\psi_4} = -2$. The formulation also sets the upper bounds of optimization variables as $U_\theta = 1$ (meaning the highest possible stiffness value is twice the initial/nominal value) and $U_{\psi_4} = 2$. All illustrations in Section 3 come from this numerical example. In particular, Eq. (6) shows the numeric expression of the model updating problem. Eq. (10) illustrates the base monomial vectors for the SOS decomposition of the objective function. Eq. (12) is the primal semidefinite programming (SDP) problem generated by SOS method from Eq. (6). Eq. (18) is the dual SDP problem generated by SOS method from Eq. (6).

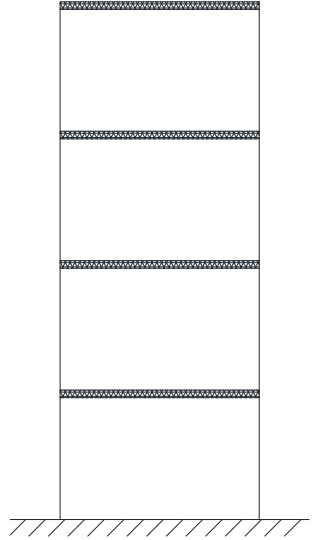

FIGURE 2 A four-story shear frame

FIGURE 3 illustrates the nonconvex objective function $f(\mathbf{x})$ against the two variables, $\theta$ and $\psi_4$. FIGURE 3(a) plots the contour of objective function over the entire feasible set $\{(\theta, \psi_4)|-1 \leq \theta \leq 1, -2 \leq \psi_4 \leq 2\}$. The global optimum $\mathbf{x}^* = (\theta^*, \psi_4^*)$ is at $(-0.100, 1.154)$, which corresponds to the "true" values of the two variables, i.e. the ideal solution. Two local optimal points, named as $\mathbf{x}_{GN} = (-1.000, 0.827)$ and $\mathbf{x}_{TR} = (-1.000, 0.000)$, locate at the boundary. This contour plot clearly shows that the objective function is nonconvex, especially around the squared region where a saddle point $s = (-0.944, 1.000)$ is. FIGURE 3 (b) shows the 3D close-up of $f(\mathbf{x})$ around the saddle point with the vertical axis as the objective function value. The figure again demonstrates the nonconvexity of this small model updating problem.

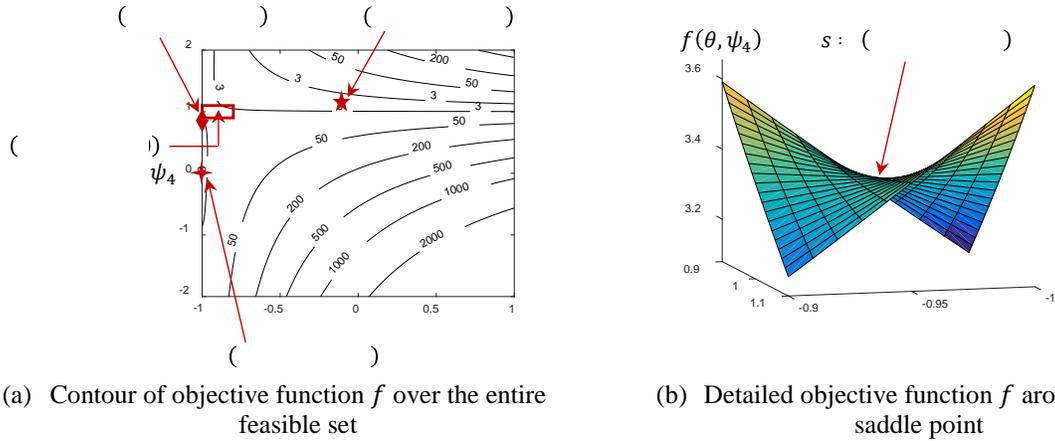

(a) Contour of objective function $f$ over the entire feasible set

(b) Detailed objective function $f$ around the saddle point

FIGURE 3 Plot of objective function $f(\mathbf{x})$, i.e. $f(\theta, \psi_4)$

By using SOS optimization method, the nonconvex problem is recast into a convex SDP problem using the formulation in Eq. (13), and the dual problem is illustrated in Eq. (16). By solving the primal and dual problems in Eq. (13) and Eq. (16), the optimal solutions can be calculated as $\gamma^* = 0.000$ for the primal and $\mathbf{y}^* = (1, -0.100, 1.154, \cdots)$ for the dual. Recalling Eq. (17), the optimal solution $\mathbf{x}^*$ for problem (6) is now easily extracted as $(-0.100, 1.154)$, which is the same as the global optimum shown in FIGURE 3.

To compare with the SOS optimization method, two local optimization algorithms are adopted to solve the optimization problem. The first local optimization algorithm is Gauss-Newton algorithm for nonlinear least squares problems [30]. Gauss-Newton algorithm is a modified version of Newton algorithm with an approximation of the Hessian matrix by omitting the higher order term. Through the MATLAB command `lsqnonlin` [31], the second algorithm uses trust region reflective algorithm [32]. The algorithm heuristically minimizes the objective function by solving a sequence of quadratic subproblems subject to ellipsoidal constraints.

For a nonconvex problem, depending on different search starting points, a local optimization algorithm may converge to different locally optimal points. TABLE 1 summarizes the model updating results calculated by different algorithms. The results show that if the search starting point happens to be close to the saddle point in FIGURE 3(b), Gauss-Newton algorithm and the trust region reflective algorithm converge at boundary points $\mathbf{x}_{\text{GN}}$ and $\mathbf{x}_{\text{TR}}$, respectively. The corresponding objective function values are both much larger than $\gamma^* = 0.000$. Only when the starting point is luckily chosen to be far away from the saddle point, the local optimization algorithms can find the global optimum. On the other hand, the SOS optimization method does not require any search starting point on $(\theta, \psi_4)$, but recasts the nonconvex problem into a convex optimization problem and reliably reaches the global optimum.

TABLE 1 Model updating results

| Optimization algorithms | Starting point | | | Updated value | | | Error | | Objective function value |
|---|---|---|---|---|---|---|---|---|---|
| | | $\theta$ | $\psi_4$ | | $\theta$ | $\psi_4$ | $\theta$ | $\psi_4$ | |
| Gauss-Newton | around $s$ | $-0.950$ | $1.000$ | $\mathbf{x}_{\text{GN}}$ | $-1.000$ | $0.827$ | $-100\%$ | $-28.30\%$ | $2.898$ |
| | $\mathbf{x}_0$ | $0.000$ | $0.000$ | $\mathbf{x}^*$ | $-0.100$ | $1.154$ | $0.00\%$ | $0.00\%$ | $0.000$ |
| Trust region reflective | around $s$ | $-0.950$ | $1.000$ | $\mathbf{x}_{\text{TR}}$ | $-1.000$ | $0.000$ | $-100\%$ | $-100\%$ | $1.914$ |
| | $\mathbf{x}_0$ | $0.000$ | $0.000$ | $\mathbf{x}^*$ | $-0.100$ | $1.154$ | $0.00\%$ | $0.00\%$ | $0.000$ |
| SOS | | – | – | $\mathbf{x}^*$ | $-0.100$ | $1.154$ | $0.00\%$ | $0.00\%$ | $0.000$ |

### 4.2 Experimental validation

To further validate the proposed model updating approach, laboratory experiments are conducted. A four-story 2D shear frame structure is used for the experiment [33]. Floor masses are accurately weighed, but all four inter-story stiffness values require updating. The structure is mounted on a shake table which provides base excitation for the structure. Accelerometers are installed on all the floors (#1~#4) and the base (#0), as shown in FIGURE 4. Modal properties of the four-story structure are then extracted from the experimental data. To conduct the structural model updating, two cases are studied. In case 1, the experimental data at all DOFs are used, which effectively makes the modal dynamic residual approach a convex optimization problem (see explanation after Eq. (3)). In addition, it is assumed all the four modes are measured/available for model updating, i.e. $n_{\text{modes}} = 4$. In case 2, the acceleration data only at the first three floors (#1~#3) are used, coinciding the instrumentation scenario in Section 4.1 – Numerical simulation. However, to make the problem scale larger as enabled by more data, it is assumed the first two modes are measured/available for model updating, i.e. $n_{\text{modes}} = 2$. Because of partial instrumentation in case 2, the model updating problem in Eq. (3) leads to a nonconvex optimization problem.

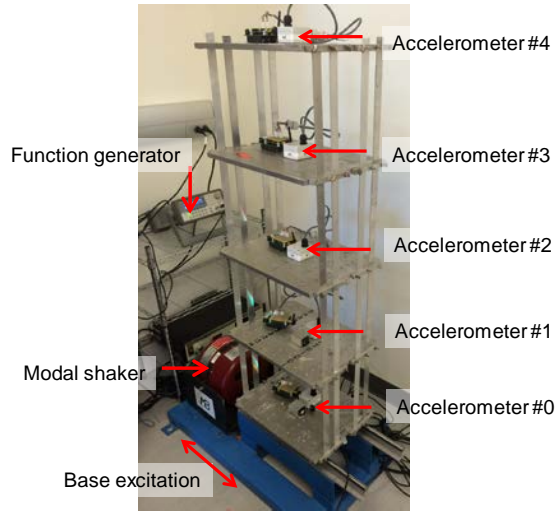

FIGURE 4 Experimental setup

### 4.2.1 Case 1: complete measurement

To obtain the modal properties of the four-story structure, a chirp signal (increasing from 0 Hz to 15 Hz in 15 seconds) is generated as base excitation. The eigensystem realization algorithm (ERA) [34] is applied to extract the modal properties from the acceleration data. FIGURE 5 shows the extracted resonance frequencies $\omega_i$ and mode shapes $\psi_i$, $i = 1,2,3,4$ for the four modes of the structure.

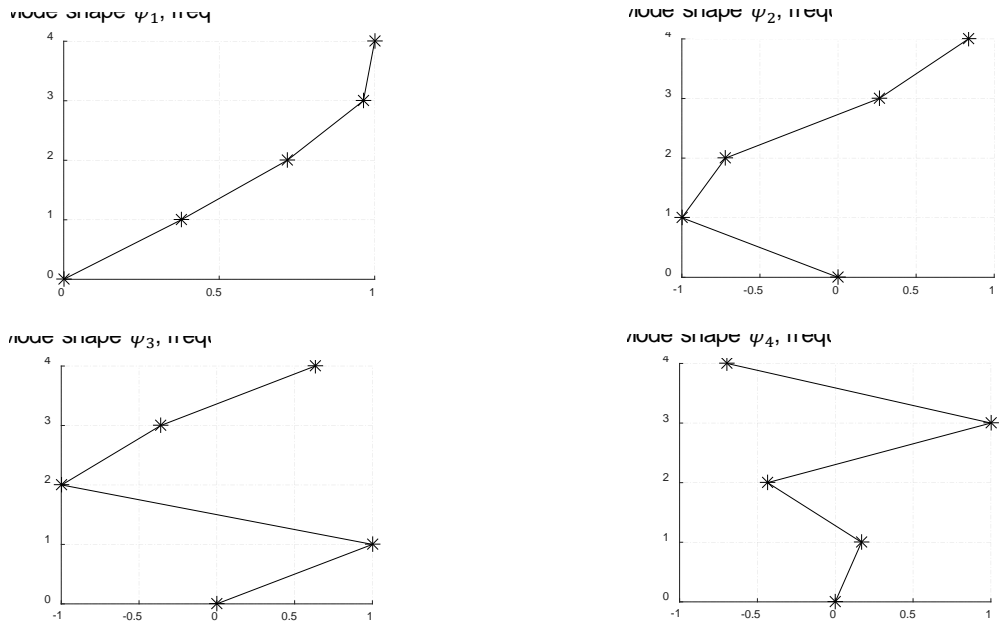

FIGURE 5 First four mode shapes extracted from chirp excitation data of the four-story structure

All the four initial/nominal inter-story stiffness values are set as 10 lbf/in prior to model updating. To update all the four stiffness values, four stiffness updating variables $\boldsymbol{\theta} \in \mathbb{R}^4$, are among the optimization variables. Same as in the simulation, each $\theta_i$ represents the relative change of a stiffness value, i.e. $\theta_i = (k_i - 10)/10$, $i = 1,2,3,4$. The weight of each floor includes the aluminum plate and the installed sensor. Considered as accurate, each floor is weighed by a scale and found to be 12.060 lb.

Because all DOFs are instrumented/measured, the modal dynamic residual approach (Eq. (3)) degenerates to a convex optimization problem as follows. Since all four modes are assumed available for experimental data, $n_{\text{modes}} = 4$.

$$\underset{\boldsymbol{\theta}}{\text{minimize}} \quad \sum_{i=1}^{n_{\text{modes}}} \left\| \left[ \mathbf{K}(\boldsymbol{\theta}) - \omega_i^2 \mathbf{M} \right] \boldsymbol{\psi}_i \right\|_2^2 \qquad (20)$$
$$\text{subject to} \quad \mathbf{L}_{\boldsymbol{\theta}} \leq \boldsymbol{\theta} \leq \mathbf{U}_{\boldsymbol{\theta}}$$

The lower bound $\mathbf{L}_{\boldsymbol{\theta}} = -\mathbf{1}_{4\times 1}$ and the higher bound $\mathbf{U}_{\boldsymbol{\theta}} = \mathbf{1}_{4\times 1}$. For comparison, the local optimization algorithms introduced in Section 4.1 are also adopted to solve the optimization problem. The starting point of $\boldsymbol{\theta}$ is selected as $\boldsymbol{\theta}_0 = (0 \quad 0 \quad 0 \quad 0)^T$. The updated inter-story stiffness values obtained from different optimization methods are summarized in TABLE 2. The updating results obtained from different optimization methods are consistent, because of the convexity of the problem in Eq. (20). Due to significant $P$-$\Delta$ effect of the lab structure, lower stories demonstrate much less inter-story stiffness.

TABLE 2 Updated inter-story stiffness values (Unit: lbf/in)

| Parameter | Gauss-Newton | Trust-region-reflective | SOS |
|---|---|---|---|
| $k_1$ | 6.949 | 6.949 | 6.949 |
| $k_2$ | 8.103 | 8.103 | 8.103 |
| $k_3$ | 9.094 | 9.094 | 9.094 |
| $k_4$ | 14.650 | 14.650 | 14.650 |

TABLE 3 compares the modal properties extracted from experimental data and generated by numerical models with the initial and the updated stiffness values. Because the updated stiffness values are the same for different optimization algorithms, only one column of "updated model" is provided in TABLE 3 for all algorithms. Besides resonance frequencies $f_n$, the modal assurance criterion (MAC) value qualifies the similarity between each experimentally extracted mode shape and the one generated by the corresponding numerical model. Essentially the square of cosine of the angle between two mode shape vectors, a MAC value closer to 1 indicates higher similarity. The natural frequencies and mode shapes of the updated model (from all the three optimization algorithms) are generally closer to the experimental results than these of the initial model. Using both the initial and the updated stiffness numbers, the values of the objective function in Eq. (20) are listed at the bottom of TABLE 3. The minimized value (10.380) of the objective function is much smaller than the value calculated using the initial stiffness numbers (271.974).

TABLE 3 Comparison of model updating results with complete measurement

| Modes | Experimental results | Initial model | | | Updated model | | |
|---|---|---|---|---|---|---|---|
| | $f_n$ (Hz) | $f_n$ (Hz) | $\Delta f_n$ (%) | MAC | $f_n$ (Hz) | $\Delta f_n$ (%) | MAC |
| 1st mode | 0.88 | 0.99 | 12.02% | 1.00 | 0.88 | 0.70% | 1.00 |
| 2nd mode | 2.75 | 2.85 | 3.64% | 0.96 | 2.74 | 0.45% | 0.98 |
| 3rd mode | 4.30 | 4.36 | 1.47% | 0.74 | 4.29 | 0.13% | 1.00 |
| 4th mode | 5.53 | 5.35 | 3.21% | 0.71 | 5.53 | 0.04% | 0.98 |
| Objective value in Eq. (20) | | 271.974 | | | 10.380 | | |

### 4.2.2 Case 2: partial measurement

In case 2, the acceleration data at only the first three floors (#1~#3) are assumed to be available and it is assumed the first two modes are measured/available for model updating, i.e. $n_{modes} = 2$, $\psi_{1,m}$ and $\psi_{2,m} \in \mathbb{R}^3$. The unmeasured/uninstrumented mode shape entry in the first mode is denoted as $\psi_{1,u} = \psi_{1,4}$ (4th-DOF); the unmeasured entry in the second mode is denoted as $\psi_{2,u} = \psi_{2,4}$ (4th-DOF). FIGURE 6 demonstrates that the resonance frequencies and the first three entries in both modes, where the known quantities are the same as those in FIGURE 5.

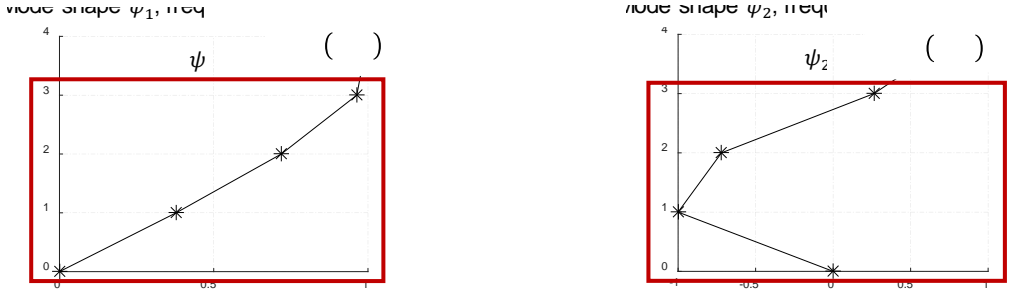

FIGURE 6 First two modes of the four-story structure; in Case 2, it's assumed only the three lower DOFs are instrumented

In this model updating problem, all the four inter-story stiffness variables $\boldsymbol{\theta} \in \mathbb{R}^4$ and the unmeasured entry in two mode shapes, $\boldsymbol{\psi}_u = (\psi_{1,4}, \psi_{2,4})^T \in \mathbb{R}^2$, are formulated as optimization variables, i.e. $\mathbf{x} = (\boldsymbol{\theta}, \boldsymbol{\psi}_u)^T \in \mathbb{R}^6$. The mass values and the initial inter-story stiffness values are the same as those in Case 1. The three optimization algorithms are again applied to solve the model updating problem, i.e. Gauss-Newton, trust region reflective, and the proposed SOS method. For the two local optimization methods, 1,000 search starting points of the updating variables $\mathbf{x} = (\boldsymbol{\theta}, \boldsymbol{\psi}_u) \in \mathbb{R}^6$ are randomly and uniformly generated between the lower bound $\mathbf{L_x} = (-1 \quad -1 \quad -1 \quad -1 \quad -2 \quad -2)^T$ and the upper bound $\mathbf{U_x} = (1 \quad 1 \quad 1 \quad 1 \quad 2 \quad 2)^T$. Starting from each of the 1,000 points, both local optimization algorithms are used to search the optimal solution. FIGURE 7 plots the optimized objective function values from 1,000 starting points by each local optimization algorithm. FIGURE 7(a) plots the performance of Gauss-Newton algorithm. The plot shows that many of the final solutions (785 out of 1,000) converge at the lowest minimum point, with the value of objective function as 0.161. However, some local optimal points are quite far away from the minimum point, and the achieved values of objective function are much higher than 0.161. FIGURE 7(b) shows the performance of trust-region-reflective algorithm. It turns out that the final solutions are separated into two groups. Many of the final solutions (733 out of 1,000) converge at the lowest minimum point with the values of objective function as 0.161. However, all the other 267 optimal points end at a local minimum with the values of objective function as 16.526.

TABLE 4 compares the modal properties obtained from experimental data and numerical models with initial and two sets of updated parameters. One updated model uses the stiffness parameters identified by SOS method. The model provides resonance frequencies and mode shapes that are very close to the experimented modal properties. The value of objective function is found to be 0.161, much smaller than that of the initial model (827.420). The other updated model uses the stiffness parameters obtained from the 267 trials by trust-region-reflective algorithm that end at the local minimum with objective function value of 16.526. Although the local minimum has a smaller objective value than that of initial model, the modal properties, including resonance frequencies and mode shapes, from the updated model are quite different from the experimental ones. This confirms that the updated parameters at a local minimum are far away from the true values, and that a local optimization algorithm cannot guarantee global optimality.

Although for a small problem like this, more randomly generated starting points can increase the chance of finding global optimum, the strategy may become challenged by larger problems, with much higher nonconvexity. On the other hand, the SOS optimization method recasts the original problem as a convex SDP and can reliably find the lowest minimum point, without searching from a large quantity of randomized starting points.

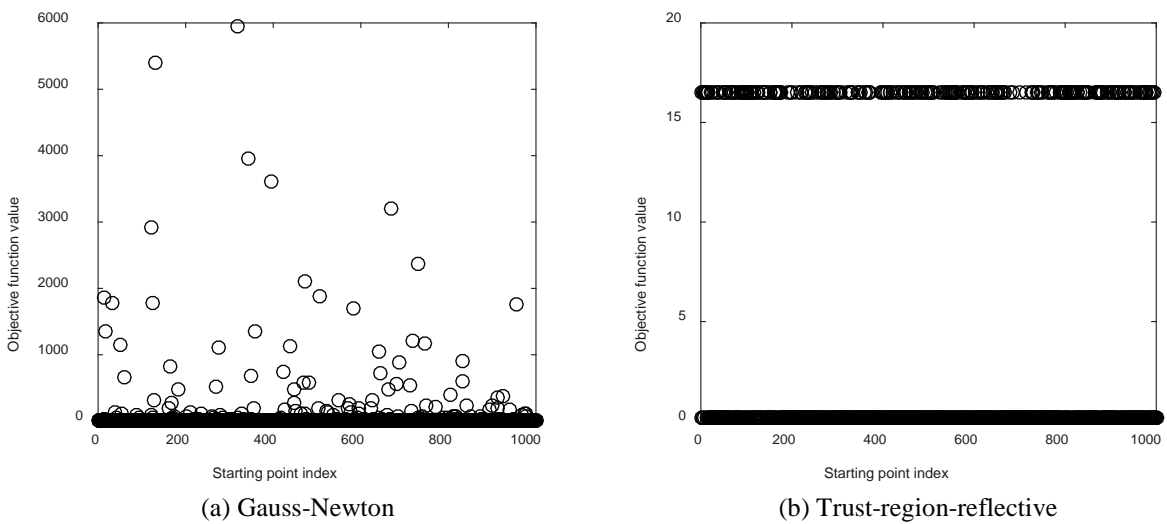

(a) Gauss-Newton  (b) Trust-region-reflective

FIGURE 7 Optimized objective function value

TABLE 4 Comparison of model updating results with partial measurement

| Modes | Experimental results | Initial model | | | SOS method | | | Local minimum from 267 trust-region-reflective trials | | |
|---|---|---|---|---|---|---|---|---|---|---|
| | $f_n$ (Hz) | $f_n$ (Hz) | $\Delta f_n$ (%) | MAC | $f_n$ (Hz) | $\Delta f_n$ (%) | MAC | $f_n$ (Hz) | $\Delta f_n$ (%) | MAC |
| 1st mode | 0.88 | 0.99 | 12.02% | 1.00 | 0.87 | 1.05% | 1.00 | 0.00 | 100% | 0.38 |
| 2nd mode | 2.75 | 2.85 | 3.64% | 0.96 | 2.75 | 0.04% | 1.00 | 0.93 | 66.03% | 0.09 |
| Objective value in Eq. (3) | | 827.420 | | | 0.161 | | | 16.526 | | |

## 5   Conclusion and discussion

This paper investigates sum of squares (SOS) optimization method for structural model updating with modal dynamic residual formulation. The modal dynamic residual formulation (Eq. (3)) has a polynomial objective function and polynomial constraints. Local optimization algorithms can be applied directly to the optimization problem, while they cannot guarantee to find the global optimum. The SOS optimization method can recast a nonconvex polynomial optimization problem into a convex semidefinite programming (SDP) problem, for which off-the-shelf solvers can reliably find the global optimum.

In particular, the nonconvex optimization problem (Eq. (3)) is recast into convex SDP problems (Eq. (13) and (16)). By solving these two SDP problems, the best lower bound and the minimizer of the original model updating problem can be solved reliably. Numerical simulation and laboratory experiments are conducted to validate the proposed model updating approach. It is shown that compared with local optimization algorithms, the proposed approach can reliably find the lowest minimum point for both complete and partial measurement cases. In addition, to improve the chance of finding a better solution, traditionally used local algorithms usually need to search from a large quantity of randomized starting points, yet the proposed SOS method does not need to.

It should be clarified that this research mainly focuses on using the SOS method to solve model updating problems towards structural system identification (instead of damage detection). In this preliminary work, the model updating of a simple shear-frame structure has been studied for a clear understanding of the proposed SOS method. Despite the advantage of reliably finding the global minimum, the SOS method is limited to optimization problems described by polynomial functions. To extend the application to other model updating formulations with non-polynomial functions, techniques such as introducing auxiliary variables and equality constraints can be utilized to convert the functions into polynomials [35]. Future work can focus on applying the SOS method on the model updating of large structures, as well as model updating formulations described by non-polynomial functions.

## ACKNOWLEDGEMENTS
This research was partially funded by the National Science Foundation (CMMI-1150700) and the China Scholarship Council (#201406260201). Any opinions, findings, and conclusions or recommendations expressed in this publication are those of the authors and do not necessarily reflect the view of the sponsors.